\documentclass[10pt,draft]{amsart}
\usepackage{amsmath,amssymb,amsrefs}

\numberwithin{equation}{section}

\newtheorem{theorem}[equation]{Theorem}
\newtheorem{lemma}[equation]{Lemma}
\newtheorem{corollary}[equation]{Corollary}
\newtheorem{prop}[equation]{Proposition}

\theoremstyle{definition}
\newtheorem{definition}[equation]{Definition}
\newtheorem{example}[equation]{Example}

\theoremstyle{remark}
\newtheorem{remark}[equation]{Remark}

\begin{document}

\begin{center}
\texttt{Comments, suggestions, corrections, and further references
  are most welcomed!}
\end{center}
\bigskip

\title{Capable groups of prime exponent and class two}
\author{Arturo Magidin}
\address{Department of Mathematical Sciences, The University of
Montana, Missoula MT 59812}
\email{magidin@member.ams.org}


\begin{abstract}
A group is called capable if it is a central factor group. We consider
the capability of finite groups of class two and exponent~$p$, $p$ an
odd prime. We restate the problem of capability as a problem about
linear transformations, which may be checked explicitly for
any specific instance of the problem. We use this restatement to
derive some known results, and prove new ones. Among them, we reduce
the general problem to an oft-considered special case, and prove that
a $3$-generated group of class~$2$ and exponent~$p$ is either cyclic
or~capable.
\end{abstract}

\maketitle

\section{Introduction}\label{sec:intro}

In his landmark paper~\cite{hallpgroups} on the classification of finite
$p$-groups, P.~Hall remarked:
\begin{quote}
The question of what conditions a group $G$ must fulfill in order that
it may be the central quotient group of another group $H$, $G\cong
H/Z(H)$, is an interesting one. But while it is easy to write down a
number of necessary conditions, it is not so easy to be sure that they
are sufficient.
\end{quote}
Following M.~Hall and Senior~\cite{hallsenior}, we make the following
definition:

\begin{definition} A group $G$ is said to be \textit{capable} if and
  only if
there exists a group $H$ such that $G\cong H/Z(H)$; equivalently, if
and only if $G$ is isomorphic to the inner automorphism group of a
group~$H$.
\end{definition}

Capability of groups was first studied by R.~Baer in~\cite{baer},
where, as a corollary of deeper investigations, he characterised the
capable groups that are a direct sum of cyclic groups.  Capability of
groups has received renewed attention in recent years, thanks to
results of Beyl, Felgner, and Schmid~\cite{beyl} characterising the
capability of a group in terms of its epicenter; and more recently to
work of Graham Ellis~\cite{ellis} that describes the epicenter in
terms of the nonabelian tensor square of the group.

We will consider here the special case of nilpotent groups of class
two and exponent an odd prime~$p$. This case was studied
in~\cite{heinnikolova}, and also addressed elsewhere (e.g., Prop.~9
in~\cite{ellis}). As noted in the final paragraphs
of~\cite{baconkappe}, currently available techniques seem insufficient
for a characterization of the capable finite $p$-groups of class~$2$, but a
characterization of the capable finite groups of class~$2$ and
exponent~$p$ seems like a more modest and possibly attainable
goal. The present work is a contribution towards achieving that goal.

We will use some of the techniques used in~\cite{capable}. We prove
that the problem of capability is equivalent to a linear algebra
problem, and then exploit this equivalence. In the rest of the
introduction we will recall the relevant definitions and basic
results. In Section~\ref{sec:candidate} we construct a natural
candidate to witness the capability of a given group~$G$. Then we show
that the properties of this natural candidate may be codified using
linear algebra in Section~\ref{sec:linear}, and deduce the promised
equivalence. In Section~\ref{sec:consequences} we exploit this
equivalence and prove several results on capability, among them
reducing the general case to the study of groups~$G$ which satisfy
$[G,G]=Z(G)$, an oft-considered special case. In
Section~\ref{sec:geometry} we present a geometric argument which
establishes some further cases of capability by using algebraic
geometry.

Throughout the paper, $p$ will be an odd prime. All groups will be
written multiplicatively, and the identity element will be denoted by
$e$; if there is danger of ambiguity or confusion, we will use $e_G$
to denote the identity of the group~$G$. The center of $G$ is denoted
by $Z(G)$. Recall that if $G$ is a group, and $x,y\in
G$, we let the commutator of $x$ and $y$ be $[x,y]=x^{-1}y^{-1}xy$. We
write commutators left-normed, so that
$[x,y,z] = [[x,y],z]$.
Given subsets $A$ and $B$ of~$G$ we define $[A,B]$ to be the subgroup
of $G$ generated by all elements of the form $[a,b]$ with $a\in A$,
$b\in B$.  The terms of the lower central series of $G$ are defined
recursively by letting $G_1=G$, and $G_{n+1}=[G_n,G]$.  A group is
\textit{nilpotent of class at most~$k$} if and only if
$G_{k+1}=\{e\}$, if and only if $G_k\subset Z(G)$. We usually drop the
``at most'' clause, it being understood. The class of all nilpotent
groups of class at most~$k$ is denoted by~$\mathfrak{N}_k$.

The following commutator identities are well known, and may be
verified by direct calculation:
\begin{prop} Let $G$ be any group. Then for all $x,y,z\in G$,
\begin{itemize}
\item[(a)] $[xy,z] = [x,z][x,z,y][y,z]$.
\item[(b)] $[x,yz] = [x,z][z,[y,x]][x,y]$.
\item[(c)] $[x,y,z][y,z,x][z,x,y] \equiv e \pmod{G_4}$.
\item[(d)] $[x^r,y^s] \equiv
  [x,y]^{rs}[x,y,x]^{s\binom{r}{2}}[x,y,y]^{r\binom{s}{2}}
  \pmod{G_4}$.
\item[(e)] $[y^r,x^s] \equiv
  [x,y]^{-rs}[x,y,x]^{-r\binom{s}{2}}[x,y,y]^{-s\binom{r}{2}}
  \pmod{G_4}$.
\end{itemize}
Here, $\binom{n}{2} = \frac{n(n-1)}{2}$ for all integers~$n$.
\label{prop:commident}
\end{prop}

As in~\cite{capable}, our main tool will be the nilpotent product of
groups, specifically the $2$-nilpotent and $3$-nilpotent products. We
restrict Golovin's original definition \cite{golovinnilprods} to the
situation we will consider:

\begin{definition} Let $A_1,\ldots,A_n$ be cyclic groups. The
  $k$-nilpotent product of $A_1,\ldots,A_n$, denoted by
  $A_1\amalg^{\germ N_k}\cdots \amalg^{\germ N_k} A_n$, is defined to
  be the group $G=F/F_{k+1}$, where $F$ is the free product of the
  $A_i$, $F=A_1*\cdots*A_n$, and $F_{k+1}$ is the $(k+1)$-st term of
  the lower central series of~$F$.
\end{definition}

Note that if $G$ is the $k$-nilpotent product of the $A_i$, then
$G\in{\germ N}_k$, and 
$G/G_k$ is the $(k-1)$-nilpotent product of the $A_i$.

\begin{theorem}[R.R.~Struik; Theorem~3 in~\cite{struikone}] Let
  $A_1\,\ldots,A_t$ be cyclic groups of order
  $\alpha_1,\ldots,\alpha_t$, respectively; if $A_i$ is infinite
  cyclic, let $\alpha_i=0$. Let $x_i$ generate $A_i$, and let $F$ be
  their free product
\[ F = A_1 * \cdots * A_t.\]
Assume that all primes appearing in the factorizations of the
$\alpha_i$ are greater than or equal to $3$. Let
$\alpha_{ij}=\gcd(\alpha_i,\alpha_j)$,
$\alpha_{ijk}=\gcd(\alpha_i,\alpha_j,\alpha_k)$ for all $1\leq i,j,k\leq
t$. Then every $g\in F/F_{4}$ can be uniquely expressed as
\begin{eqnarray*}
 g &=& x_1^{a_1}\cdots x_t^{a_t}\prod_{1\leq i<j\leq
   t}[x_j,x_i]^{a_{ji}}\prod_{1\leq i<j\leq
   t}[x_j,x_i,x_i]^{a_{jii}}[x_j,x_i,x_j]^{a_{jij}}\\
&&\qquad\qquad\qquad\prod_{1\leq i<j<k\leq
   t}[x_j,x_i,x_k]^{a_{jik}}[x_k,x_i,x_j]^{a_{kij}}
\end{eqnarray*}
where the $a_i$ are taken modulo $\alpha_i$, the $a_{ji}$ modulo
$\alpha_{ji}$, and the $a_{{\ell}mn}$ modulo $\alpha_{{\ell}mn}$.
\end{theorem}

A collection process may be used to obtain a multiplication table for
this group, although this will not be necessary for our purposes. We
direct the reader to Struik's original paper~\cite{struikone} (where
she uses a slightly different choice of normal form), or
to~\cite{capable} (where the multiplication formulas agree
with the normal form given above).

To obtain a description of the $2$-nilpotent product, already given by
Golovin~\cite{metab}, we simply take the quotient modulo the third
term of the lower central series; this is the subgroup of all elements
of the form
\[\prod_{1\leq i<j\leq
   t}[x_j,x_i,x_i]^{a_{jii}}[x_j,x_i,x_j]^{a_{jij}} \prod_{1\leq i<j<k\leq
   t}[x_j,x_i,x_k]^{a_{jik}}[x_k,x_i,x_j]^{a_{kij}}\]
for integers $a_{jii}$, $a_{jij}$, $a_{kij}$, and $a_{jik}$.

From the definition, it is clear that the $k$-nilpotent product is the
coproduct in the variety~${\germ N}_k$, so it will have the usual universal
property. When we take the $k$-nilpotent product of infinite cyclic groups
we obtain the relatively free group in the variety ${\germ
  N}_k$.

Finally, when we say that a group is \textit{$k$-generated} we mean that it
can be generated by $k$ elements, but may in fact need less. If we
want to say that it can be generated by $k$ elements, but not by $m$ elements
for some $m<k$, we will say that it is \textit{minimally
$k$-generated}, or \textit{minimally generated by $k$ elements}.

\section{A natural candidate}\label{sec:candidate}

Using the $2$- and $3$-nilpotent product of cyclic groups, we can
produce a natural ``candidate for witness'' to the capability of a
given finite nilpotent groups of class~$2$ and exponent~$p$, which we
can then investigate directly.  This will be the theme throughout this
article.  We begin with an easy observation:

\begin{lemma}[cf.~Lemma~2.1 in~\cite{isaacs}]
Let $G$ be a capable group, generated by $g_1,\ldots,g_n$, $n>1$. Then
there exists a group $H$, and elements $h_1,\ldots,h_n\in H$ such that
$H$ is generated by $h_1,\ldots,h_n$, and $H/Z(H)\cong G$, where the
isomorphism is induced by the map sending $h_i$ to the corresponding $g_i$.
\end{lemma}
\begin{proof} Let $K$ be any group such that $K/Z(K)\cong G$. Let
  $h_i\in K$ be any element mapping to $g_i$. Let $H=\langle
  h_1,\ldots,h_n\rangle$. We then have $HZ(K)=K$, so it follows that
  $Z(H)=Z(K)\cap H$ and hence $H/Z(H)\cong H/(Z(K)\cap H) \cong K/Z(K) \cong G$.
\end{proof}

Let $G$ be a finite nilpotent group of class at most~$2$ and
exponent~$p$, minimally generated by $g_1,\ldots,g_n$, $n>1$. Then $G$
is a quotient of the $2$-nilpotent product of $n$ cyclic $p$-groups,
\[ G \cong \left(\langle x_1\rangle \amalg^{{\germ N}_2}\cdots
\amalg^{{\germ N}_2}\langle x_n\rangle\right)/N\] where $N$ is the
kernel of the map induced by mapping $x_i$ to $g_i$ (using the
universal property of the coproduct), and is contained
in the commutator subgroup of the $2$-nilpotent product. Let
$y_1,\ldots,y_n$ generate infinite cyclic groups, let $\germ K$ be the
$3$-nilpotent product of the cyclic groups,
\[{\germ K} = \langle y_1\rangle \amalg^{{\germ N}_3}\cdots
\amalg^{{\germ N}_3} \langle y_n\rangle,\]
and let $\mathcal{K} = {\germ K}/[{\germ K},{\germ K}]^p$.  By abuse
of notation, we denote the images of the $y_i$ in $\mathcal{K}$ by
$y_i$ as well.

Note that the commutator subgroup of $\langle x_1\rangle
\amalg^{{\germ N}_2}\cdots \amalg^{{\germ N}_2} \langle x_r\rangle$ is
isomorphic to the subgroup of $[\mathcal{K},\mathcal{K}]$ generated by
the commutators of the form $[y_j,y_i]$, $1\leq i<j\leq n$, by mapping
$[x_j,x_i]$ to $[y_j,y_i]$. Let $\mathcal{N}$ be the subgroup of
$\mathcal{K}$ that corresponds to $N$ under this identification.

Let $\mathcal{M}=[\mathcal{N},\mathcal{K}]$. Finally, let
$K = \mathcal{K}/\mathcal{M}$. Again, we also denote the images of the
$y_i$ in $K$ by $y_i$.

\begin{theorem}
Let $G$, ${\germ
    K}$, $\mathcal{K}$, and~$K$ be as in the previous three 
  paragraphs. Then $G$ is capable if and only if $G\cong K/Z(K)$.
\end{theorem}

\begin{proof} First, note that
  $K/\langle K^p,\mathcal{N}\mathcal{M},K_3\rangle$ is isomorphic to~$G$.
  Since each of  $K^p$, $\mathcal{N}\mathcal{M}$, and~$K_3$ are
  central in~$K$, it is
  always be the case that $K/Z(K)$ is a quotient of~$G$.

Assume $G$ is capable, and let $H \in {\germ N}_3$ be a
group with $H/Z(H)\cong G$. We may assume $H$ is generated by
elements $h_1,\ldots,h_n$ mapping to $g_1,\ldots,g_n$,
respectively. Then $H^p\subset Z(H)$. We claim this implies that
$H_2^p=\{e\}$.

Indeed, first note that $H_3^p=\{e\}$. For if $c\in H_2$, $h\in H$,
then using Proposition~\ref{prop:commident} we have that $[c,h]^p = [c,h^p]
= e$. Since $H_3$ is abelian, and generated by all elements $[c,h]$ as
above, this shows that $H_3$ is of exponent~$p$. Now let $a,b\in H$; we have (once again using
Proposition~\ref{prop:commident}):
\begin{eqnarray*}
\relax [a,b]^p & = &
       [a^p,b][a,b,a]^{-\binom{p}{2}}[a,b,b]^{-p\binom{1}{2}}\\
& = & [a^p,b][a,b,a]^{-p\left(\frac{p-1}{2}\right)}\\
& = & e.
\end{eqnarray*}
Since $H_2$ is abelian and generated by all such $[a,b]$, which are of
exponent~$p$, we conclude that $H_2^p=\{e\}$, as claimed.

The natural surjection ${\germ K}\to G$, mapping $y_i$ to $g_i$
therefore factors through $H$ and~$\mathcal{K}$, so we have the following
exact diagram:
\[
\begin{array}{ccccccccc}
&&&&\mathcal{K}\\
&&&&\downarrow\\
1&\rightarrow& Z(H)&\rightarrow & H & \rightarrow& G & \rightarrow &
1.\\
&&&&\downarrow\\
&&&&1
\end{array}\]
If $y\in\mathcal{N}$, then its image in $H$ must map to the trivial
element in~$G$, so its image is in the center of $H$. It follows that
$\mathcal{M}= [\mathcal{N},\mathcal{K}]$ is in the kernel of the map from
$\mathcal{K}$ to $H$, so the map factors through $K$ giving a diagram
with exact row and column:
\[
\begin{array}{ccccccccc}
&&&&K\\
&&&&\downarrow\\
1&\rightarrow& Z(H)&\rightarrow & H & \rightarrow& G & \rightarrow & 1.\\
&&&&\downarrow\\
&&&&1\\
\end{array}\]
Since the center of $K$ maps into the center of $H$, it follows that
$K/Z(K)$ has $G$ as a quotient. Since we already know that $K/Z(K)$ is
a quotient of $G$, it follows that $K/Z(K)\cong G$, as desired.
\end{proof}

Next, note that $Z(K/K^p) = Z(K)/K^p$, so we have:

\begin{theorem}[cf.~Theorem~8.6 in~\cite{capable}]  
Let $G$ be a nilpotent group of class at most~$2$ and exponent an odd
  prime~$p$, minimally generated by $g_1,\ldots,g_n$, with $n>1$. Let
  $A_1,\ldots,A_n$ be cyclic groups of order~$p$, generated by
  $x_1,\ldots,x_n$, respectively, and let $\mathcal{G}$ be the
  $2$-nilpotent product of the $A_i$,
$\mathcal{G} = A_1 \amalg^{{\germ N}_2}\cdots \amalg^{{\germ N}_2}
  A_n$.
Let $N\triangleleft \mathcal{G}$ be the
kernel of the homomorphism $\mathcal{G}\to G$ induced by mapping $x_i$
to $g_i$. Let $K$ be the $3$-nilpotent product of the $A_i$,
$K = A_1\amalg^{{\germ N}_3} \cdots \amalg^{{\germ N}_3} A_n$,
and identify $N$ with the corresponding subgroup of~$K_2$. Then $G$ is
capable if and only if 
\[G \cong (K/[N,K])/Z(K/[N,K]).\]
\label{thm:natcand}
\end{theorem}

What are the elements of $Z(K/[N,K])$? An element $k[N,K]$ will be
central in $K/[N,K]$ if and only if $[k,K]\subseteq [N,K]$. This
includes the center of~$K$ (which is equal to $K_3$, see Theorem~5.1
in~\cite{capable}), as well as $N$. Since $K/\langle N,K_3\rangle \cong
G$, it follows that $G$ is capable if and only if
\[ Z(K/[N,K]) = \langle K_3, N\rangle / [N,K].\]
The left hand side always contains the right hand side. However,
the problem is that there could be more elements in $Z(K/[N,K])$.
Let $k\in K$, and write $k$ in normal form:
\begin{eqnarray*}
 k &=& x_1^{a_1}\cdots x_n^{a_n}\prod_{1\leq i<j\leq
   n}[x_j,x_i]^{a_{ji}}\prod_{1\leq i<j\leq
   n}[x_j,x_i,x_i]^{a_{jii}}[x_j,x_i,x_j]^{a_{jij}}\\
&&\qquad\qquad\qquad\prod_{1\leq i<j<k\leq
   n}[x_j,x_i,x_k]^{a_{jik}}[x_k,x_i,x_j]^{a_{kij}}.
\end{eqnarray*}
Suppose that there exists $i\in \{1,\ldots,n\}$ with $a_i\neq 0$. If
 $i<n$, consider $[x_n,k]$. Using the formulas in
Proposition~\ref{prop:commident}, it is easy to verify that the normal form for
$[x_n,k]$ will have a nonzero exponent for $[x_n,x_i]$; since
 $N\subset K_2$, $[N,K]\subset K_3$, so $[x_n,k]\not\in [N,K]$. If
 $i=n$, then consider $[k,x_{n-1}]$, and once again we find that this
 element has nonzero exponent in $[x_n,x_{n-1}]$, so
 $[k,x_{n-1}]\notin [N,K]$.
That is, we have shown that if $k\in K$ satisfies
 $[k,K]\subseteq [N,K]$, then $k\in K_2$. With this observation, we
 obtain the following:

\begin{corollary}
Let $G$, $K$, and $N$ be as in Theorem~\ref{thm:natcand}. Then $G$ is
capable if and only if the following condition holds: For every
element $k\in K$ of the form
\begin{equation}
k = \prod_{1\leq i<j\leq n} [x_j,x_i]^{a_{ji}}\label{eq:kform}
\end{equation}
with $0\leq a_{ji} < p$, if $[k,x_i]\in [N,K]$ for $i=1,\ldots,n$, then
$k\in N$. 
\label{cor:candidate}
\end{corollary}

\begin{proof} We know that $G$ is capable if and only if
\[ \bigl\{ k\in K_2\,\bigm|\, [k,K]\subseteq [N,K]\bigr\} = \bigl\langle
N,K_3\rangle.\] If $k\in K_2$, then $[k,K] = \bigl\langle
[k,x_1],\ldots,[k,x_n]\bigr\rangle$ by
Proposition~\ref{prop:commident}, so if $G$ is capable then all $k$
such that $[k,x_i]\in [N,K]$, $i=1,\ldots,n$, must satisfy $k\in
\langle N,K_3\rangle$.  Conversely, assume the condition holds and let
$k\in K_2$ be such that $[k,K]\subset [N,K]$; we want to show that
$k\in\langle N,K_3\rangle$.  Multiplying by suitable elements of
$K_3$, we may assume that $k$ is of the form (\ref{eq:kform}). Since
$[k,K]\subset [N,K]$, it follows that $[k,x_i]\in [N,K]$ for each $i$;
so $k\in N$, by the condition given. This proves the corollary.
\end{proof}

\section{Some Linear Algebra}\label{sec:linear}

The main advantage of Corollary~\ref{cor:candidate}, beyond giving a
very precise and explicit condition to check, is that the condition
can also be recast as a statement about vector spaces and linear
transformations. At that point, we have a whole array of tools that
can be brought to bear upon the problem. For example,
Theorem~\ref{th:mckinnon} below uses algebraic geometry to settle
several cases.

We will now translate the condition into the promised statement about
vector spaces and linear transformations. The key observation is the
simple fact that if $H\in{\germ N}_3$, then for each $h\in H$ the map 
$\varphi_h\colon H_2\to H_3$ given by
$\varphi_h(c) = [c,h]$ is an abelian group homomorphism; if $H_2$ is
of exponent~$p$, then the map becomes a linear transformation of
spaces over $\mathbb{F}_p$, the Galois field of $p$ elements.

Let $n>1$ be a fixed integer. Let 
$K = \langle x_1\rangle \amalg^{{\germ N}_3}\cdots \amalg^{{\germ
    N}_3} \langle x_n\rangle$,
where $x_i$ is of order $p$, $p$ an odd prime. Then $K_2$ is a vector
space over $\mathbb{F}_p$, with basis given by all elements of the
form $[x_j,x_i]$ and $[x_s,x_r,x_t]$, where $1\leq i<j\leq n$, 
$1\leq r<s\leq n$, and $1\leq r\leq t\leq n$.

Let $V$ be the vector space over $\mathbb{F}_p$ of dimension
$\binom{n}{2}$, with basis 
\begin{equation}
\bigl\{v_{ji} \,\bigm|\, 1\leq i<j\leq n\bigr\}.
\label{eq:defofV}
\end{equation}
Let $W$ be a
vector space over $\mathbb{F}_p$ of dimension
$2\left(\binom{n}{2}+\binom{n}{3}\right)$ with basis
\begin{equation}
\bigl\{ w_{jik}\,\bigm|\, 1\leq i<j\leq n,\ 1\leq i\leq k\leq
n\bigr\}.
\label{eq:defofW}
\end{equation}
We will refer to these bases as the \textit{``standard bases''} for
$V$ and~$W$, and we talk about \textit{$V$ and~$W$ corresponding
to~$n$} to mean the case where the indices of the coordinates range
from $1$ to~$n$ (subject to the conditions listed above).  We can
clearly identify $K_3$ with $W$, and $K_2$ with $V\oplus W$. To
restate Corollary~\ref{cor:candidate} in terms of $V$ and $W$, we just
need to describe what happens when we take the commutator of
$[x_j,x_i]$ with~$x_r$. We have two cases: if $i\leq r$, then
$[x_j,x_i,x_r]$ is already in normal form, and there is nothing to
do. If $i>r$, then the commutator $[x_j,x_i,x_r]$ must be rewritten in
normal form. We use Proposition~\ref{prop:commident}(c): since
$[x_j,x_i,x_r][x_i,x_r,x_j][x_r,x_j,x_i] = e$, we have:
\begin{eqnarray*}
\relax[x_j,x_i,x_r] & = & [x_r,x_j,x_i]^{-1}[x_i,x_r,x_j]^{-1}\\
& = & [[x_r,x_j]^{-1},x_i][x_i,x_r,x_j]^{-1}\\
& = & [x_j,x_r,x_i][x_i,x_r,x_j]^{-1}.
\end{eqnarray*}

So, for each $r=1,\ldots,n$ we define a linear transformation
$\varphi_r\colon V\to W$ to be given by:
\begin{equation}
\varphi_r(v_{ji}) = \left\{\begin{array}{ll}
w_{jir}&\mbox{if $r\geq i$,}\\
w_{jri}-w_{irj}&\mbox{if $r<i$.}
\end{array}\right.
\label{eq:defvarphi}
\end{equation}
That is, $\varphi_r$ codifies the map $c\mapsto [c,x_r]$.  It is easy
to verify that each $\varphi_r$ is injective.

Let $G$ be a non-cyclic group of class~$2$ and exponent~$p$.  Consider
the situation in Corollary~\ref{cor:candidate}, and let $X$ be the
subspace of $V$ corresponding to the subgroup~$N$ determined by~$G$
(we will refer to such $X$ as \textit{the subspace of~$V$
corresponding to $G$} or \textit{determined by $G$}). We define $Y_X$ to
be the subspace of $W$ spanned by the images of $X$; that is:
\begin{equation}
Y_X = \Bigl\langle \varphi_1(X),\ldots,\varphi_n(X)\Bigr\rangle.
\label{eq:YsubX}
\end{equation}
Thus, $Y_X$ corresponds to the subgroup $[N,K]$. Finally, let $Z_X$ be
the subspace of~$V$ given by:
\begin{equation}
Z_X = \bigcap_{i=1}^n \varphi_i^{-1}\left(Y_X\right).
\label{eq:ZsubX}
\end{equation}
Clearly, $X\subseteq Z_X$, and $Z_X$ corresponds to the subgroup of
all $k\in K_2$, written as in (\ref{eq:kform}), such that $[k,x_i]\in
[N,K]$ for $i=1,\ldots,n$. We know that $G$ is capable if and only if
the subgroup of such $k$ equals~$N$, so we have:

\begin{theorem} Let $G$ be a finite nilpotent group of class~$2$ and
  exponent~$p$, minimally generated by $g_1,\ldots,g_n$, $n>1$. Let
  $V$ and $W$ be the vector spaces over $\mathbb{F}_p$ defined in
  $(\ref{eq:defofV})$ and~$(\ref{eq:defofW})$, let
  $\varphi_1,\ldots,\varphi_n$ be as in~$(\ref{eq:defvarphi})$, and let
  $X$ be the subspace of $V$ determined by~$G$ (that is, corresponding
  to the kernel of the natural map $\langle x_1\rangle \amalg^{{\germ
  N}_2}\cdots \amalg^{{\germ N}_2} \langle x_n\rangle \longrightarrow
  G$ given by $x_i\mapsto g_i$, where each $x_i$ is of order~$p$). Let
  $Y_X$ and $Z_X$ be defined by $(\ref{eq:YsubX})$
  and~$(\ref{eq:ZsubX})$. Then $G$ is capable if and only if $X=Z_X$.
\label{th:linalgcond}
\end{theorem}

\begin{remark}We will feel free to drop the subscript $X$ from $Y_X$ and
  $Z_X$ when it is clear from context; also, to avoid multiple
  subindices, if we have subspaces $X_1,\ldots,X_r$, we will denote
  $Y_{X_i}$ and $Z_{X_i}$ simply by $Y_i$ and $Z_i$, respectively.
\end{remark}

\section{Some consequences}\label{sec:consequences}

In this section, we will prove several results regarding the
capability of groups of class~$2$ and exponent~$p$, based on the
restatement of the problem given in Theorem~\ref{th:linalgcond}. It
will be clear that many of the results could be proven by appealing
directly to the normal forms in the groups in question, without having
to refer to linear algebra, but at least the author found that the
linear algebra setting was usually easier to think about (and see also
Section~\ref{sec:geometry}).

We start with an example of how we can use the result to prove that a
given group is not capable, and which also illustrates how to set up
the linear algebra problem given a specific group.

\begin{example} A group of class two and exponent~$p$, which is not
  capable.

Let $G=\langle x_1,x_2,x_3,x_4\rangle$ be the nilpotent group of
class two presented by: 
\begin{equation}
G = \left\langle{x_1,x_2,x_3,x_4 \left|\begin{array}{rcl}
\relax[x_3,x_1][x_3,x_2]^{-1} = [x_3,x_1][x_4,x_1]^{-1} & = & e,\\
\relax[x_4,x_2] = [x_4,x_3] = [x_2,x_1] & = & e,\\
x_1^p = x_2^p = x_3^p = x_4^p & = & e,
\end{array}\right.}\right\rangle
\label{eq:extraspecial}
\end{equation}
(we assume the condition $[G,G]\subset Z(G)$ is
given, so this is a presentation as an $\mathfrak{N}_2$-group).
This is an extra-special group of order $p^5$; it is mentioned in
Section~8 of~\cite{capable}, where the fact that it is not capable is
deduced by invoking a theorem of Beyl, Felgner, and Schmid
in~\cite{beyl}. We prove that fact here using our set-up.

Let $V$ and $W$ be the vector spaces defined in $(\ref{eq:defofV})$
and~$(\ref{eq:defofW})$,  corresponding to $n=4$; then $X$
is defined by looking at the identities defining~$G$, that is:
\[ X = \bigl\langle v_{31}-v_{32}, v_{31}-v_{41}, v_{42}, v_{43}, v_{21}
\bigr\rangle.\]
We have that ${\rm dim}(V) = \binom{4}{2} = 6$, and ${\rm
  dim}(X)=5$. So $X\neq Z_X$ if and only if $Z_X=V$. And for that, it
is enough to show that $v_{41}\in Z_X$.

To show that $v_{41}\in Z_X$, we need to show that each of $v_{411}$,
$v_{412}$, $v_{413}$, and $v_{414}$ are in 
$Y_X = \bigl\langle
\varphi_1(X),\varphi_2(X),\varphi_3(X),\varphi_4(X)\bigr\rangle$.
Indeed:
\begin{eqnarray*}
v_{411} & = & \varphi_1\bigl( v_{42} + (v_{31}-v_{32}) - (v_{31}-v_{41}) \bigr)\\
&&\qquad +\varphi_2\bigl(v_{31}-v_{41}\bigr) - \varphi_3\bigl(v_{21}\bigr) + \varphi_4\bigl(v_{21}\bigr).\\
v_{412} & = & \varphi_1\bigl(v_{42}\bigr) +
\varphi_4\bigl(v_{21}\bigr).\\
v_{413} & = & \varphi_1\bigl(v_{43}\bigr) -
\varphi_2\bigl(v_{43}\bigr) + \varphi_3\bigl(v_{42}\bigr) +
\varphi_4\bigl(v_{31}-v_{32}\bigr).\\
v_{414} & = & \varphi_3\bigl(v_{42}\bigr) - \varphi_2\bigl(v_{43}\bigr) +
\varphi_4\bigl( (v_{31}-v_{32}) -
(v_{31}-v_{41})\bigr).
\end{eqnarray*}

Therefore, $G$ is not capable. \hfill$\Box$
\end{example}

The extreme cases are easy to handle, since the $\varphi_i$ are
injective, and their images span~$W$:

\begin{lemma} Let $n>1$, and let $V$, $W$, and
  $\varphi_1,\ldots,\varphi_n$ be given as
  in~$(\ref{eq:defofV})$--$(\ref{eq:defvarphi})$. If
  $X=\{\mathbf{0}\}$, then $Z_X=X$; if $X=V$, then $Z_X=X$.
\end{lemma}

\begin{theorem} Let $n>1$, $p$ an odd prime. If $G$ is a direct sum of
  $n$ cyclic groups of order~$p$, then $G$ is capable. If $G$ is the
  $2$-nilpotent product of $n$ cyclic groups of order~$p$, then $G$ is
  capable.
\label{th:trivialcases}
\end{theorem}

\begin{proof} The direct sum corresponds to the case $X=V$, while the
  $2$-nilpotent product to $X=\{\mathbf{0}\}$. Of course, these two
  conclusions are simply  special cases of Baer's
  theorem for abelian groups and its generalization to the
  $k$-nilpotent product of cyclic $p$-groups, $k<p$ (Theorem~6.4
  in~\cite{capable}).
\end{proof}

In discussing the situation in the linear algebra setting, we will
refer to \textit{``the $(j,i)$ coordinate''} of vectors of $V$, or to
\textit{``the $(j,i,k)$ coordinate''} of vectors of~$W$. We refer, of
course, to the coefficient of $v_{ji}$ (resp.~of $w_{jik}$) when the
vector is expressed in terms of the standard bases given.

The following easy observations will be repeated several times, so we
state them as a lemma:

\begin{lemma}
Let $n>1$, and let $V$, $W$, and $\varphi_1,\ldots,\varphi_n$ be as
in~$(\ref{eq:defofV})$--$(\ref{eq:defvarphi})$.
\begin{itemize}
\item[(a)] $\varphi_r(\mathbf{v})$ has nonzero $(j,i,i)$ coordinate,
  $1\leq i<j\leq n$, if and only if $\mathbf{v}$ has nonzero $(j,i)$
  coordinate, and $r=i$.
\item[(b)] $\varphi_r(\mathbf{v})$ has nonzero $(j,i,j)$
  coordinate, $1\leq i<j\leq n$, if and only if $\mathbf{v}$ has
  nonzero $(j,i)$ coordinate, and $r=j$.
\item[(c)] $\varphi_r(\mathbf{v})$ has nonzero $(k,i,j)$ coordinate,
 $1\leq i<j\leq n$, $i<k\leq n$, if and only if either $r=j$ and
 $\mathbf{v}$ has nonzero $(k,i)$ coordinate, or else $r=i$, and
 either $k>j$ and $\mathbf{v}$ has nonzero $(k,j)$ coordinate or $j>k$
 and $\mathbf{v}$ has nonzero $(j,k)$ coordinate. In the case where
 $r=i$, the $(k,i,j)$ coordinate of $\varphi_i(\mathbf{v})$ is equal
 to minus the $(j,i,k)$ coordinate.
\end{itemize}
\end{lemma}

From these, we obtain:

\begin{lemma}
Let $n>1$, and let $V$, $W$, and $\varphi_1,\ldots,\varphi_n\colon
V\to W$ be as in~$(\ref{eq:defofV})$--$(\ref{eq:defvarphi})$.
Let $X$ be a subspace of $V$, and let $i$ and $j$ be
fixed integers, $1\leq i<j\leq n$. If all vectors in $X$ have zero
$(j,i)$ coordinate, then all vectors in $Z_X$ have zero $(j,i)$
coordinate. 
\label{lemma:zerocoord}
\end{lemma}
\begin{proof} Let $\mathbf{v}\in V$ be a vector with nonzero $(j,i)$
  coordinate. It suffices to show that $\mathbf{v}\notin
  \varphi_i^{-1}(Y_X)$. Indeed, since all vectors in $X$ have zero $(j,i)$ coordinate,
  it follows that all vectors in $Y_X$ have zero $(j,i,i)$ coordinate
  (applying part (a) of the lemma above). Since
  $\varphi_i(\mathbf{v})$ has nonzero $(j,i,i)$ coordinate, it cannot
  lie in $Y_X$, and therefore $\mathbf{v}\notin\varphi_i^{-1}(Y_X)$.
\end{proof}

We will state most of our results as lemmas about the linear algebra
situation, and deduce the corresponding conclusions for
groups of class two and exponent~$p$ as theorems. Since some of our lemmas will
appear contrived without the benefit of knowing the group theory
result we are after, we will usually state the theorem first, and then
the corresponding lemma about linear algebra.

\subsection*{A theorem of G.~Ellis}

The original impetus behind this work was my
desire to obtain an alternative proof of the following result:

\begin{theorem}[G.~Ellis, Proposition~9 of~\cite{ellis}] Let $G$ be a
  finitely generated group of nilpotency class exactly two and of
  prime exponent. Let $\{x_1,\ldots,x_k\}$ be a subset of~$G$
  corresponding to a basis of the vector space $G/Z(G)$, and
  suppose that those non-trivial commutators of the form $[x_j,x_i]$,
  $1\leq i<j\leq k$ are distinct and constitute a basis for the vector
  space $[G,G]$. Then $G$ is capable.
\label{th:ellis}
\end{theorem}

\begin{proof}
Let
$x_{k+1},\ldots,x_n$ be elements of $Z(G)$ that, together with
$x_1,\ldots,x_k$, constitute a minimal generating set of~$G$. The
conditions given in the statement of the theorem imply that the
subspace $X$ of $V$ corresponding to $G$ will have a basis consisting
of a subset of the standard basis of~$V$; therefore, Ellis's Theorem will
follow from Lemma~\ref{lemma:coordsub} below.
\end{proof}

\begin{lemma} Let $n>1$, and let $V$, $W$,
  $\varphi_1,\ldots,\varphi_n$ be given as
  in~$(\ref{eq:defofV})$--$(\ref{eq:defvarphi})$.  If $X$ has a basis
  consisting of a subset of $\{v_{ji}\,|\, 1\leq i<j\leq n\}$ (that
  is, $X$ is a ``coordinate subspace'' of~$V$), then $Z_X=X$.
\label{lemma:coordsub}
\end{lemma}

\begin{proof} Assume that $X$
  has a basis consisting of a subset $S\subset \{v_{ji}\,|\, 1
  \leq i<j\leq n\}$. If $v_{rs}\notin S$, then all vectors in $X$ have
  zero $(r,s)$ coordinate, and therefore by
  Lemma~\ref{lemma:zerocoord}, so do all vectors in $Z_X$. This
  means that $Z_X\subset \langle S\rangle = X$, giving equality.
\end{proof}

\subsection*{Capability of coproducts}

Recall that we define the $2$-nilpotent product of two ${\germ N}_2$
groups (not necessarily cyclic) $A$ and $B$ as
\[ A \amalg^{\mathfrak{N}_2} B = (A*B)/(A*B)_2,\]
where $A*B$ is their free product. A theorem of T.~MacHenry
in~\cite{machenry} shows that any element of
$A\amalg^{\mathfrak{N}_2}B$ can be written uniquely as $abc$, where
$a\in A$, $b\in B$, and $c\in [B,A]$, the `cartesian', and that the
cartesian is isomorphic to $B^{\rm ab}\otimes A^{\rm ab}$, via the
map sending $[b,a]$ to $\overline{b}\otimes \overline{a}$.

\begin{lemma} Let $G$ be a finite nontrivial group of class two and exponent
 $p$, $p$ an odd prime, and let $C_p$ be the cyclic group of order $p$.
 Then $G\amalg^{\mathfrak{N}_2} C_p$ is capable.
\label{lemma:coprodcyclic}
\end{lemma}

\begin{proof} The result follows from Theorem~\ref{th:trivialcases} if
 $G$ is cyclic, so assume that $G$ is minimally generated by
 $g_1,\ldots,g_n$, $n>1$, and denote the generator of $C_p$ by
 $x_{n+1}$. Let $V_n$ be the vector space of dimension $\binom{n}{2}$
 with basis $v_{ji}$, $1\leq i<j\leq n$, and let $V_{n+1}$ be the
 larger vector space that also includes the basis vectors
 $v_{(n+1)k}$, $1\leq k\leq n$. Identify $V_n$ with the obvious
 subspace of~$V_{n+1}$. If $X$ is the subspace of $V_n$
 corresponding to $G$, then $X$ is also the subspace of $V_{n+1}$
 corresponding to $G\amalg^{\mathfrak{N}_2} C_p$. This subspace does
 not contain any vectors with a nonzero $(n+1,k)$ coordinate, $1\leq
 k\leq n$. The result will now follow from
 Lemma~\ref{lemma:missingcoord} below.
\end{proof}

\begin{lemma} Let $n>1$, let $V$, $W$, $\varphi_1,\ldots,\varphi_n$ be
 as in $(\ref{eq:defofV})$--$(\ref{eq:defvarphi})$, and let $X$ be a
 subspace of~$V$.  Suppose that there exists $j$, $1\leq j\leq n$,
 such that for all $\mathbf{v}\in X$ and all $i$, $1\leq i\leq n$,
 $i\neq j$, $\mathbf{v}$ has zero $(i,j)$ or zero $(j,i)$ coordinate
 (whichever makes sense). Then $X=Z_X$.
\label{lemma:missingcoord}
\end{lemma}

\begin{proof} We claim that
\[ \varphi_j(V)\cap \Bigl\langle
 \varphi_1(X),\ldots,\varphi_{j-1}(X),\varphi_{j+1}(X),\ldots,\varphi_n(X)\Bigr\rangle
 = \bigl\{ \mathbf{0}\bigr\}.\]
Indeed, every nonzero coordinate of a vector in $\varphi_j(V)$ 
 has at least one index (either the second or the third) equal
 to $j$. But no vector in any of $\varphi_i(X)$, $i\neq j$, has a
 nonzero coordinate involving a $j$ (in any position). Thus, neither
 do any of their linear combinations.
This means that $\varphi_j^{-1}(Y_X) = X$ (since $\varphi_j$ is
injective), and therefore that
$Z_X = X$.
\end{proof}

In fact, we can strengthen this result considerably:

\begin{theorem} Let $G$ and $H$ be any two nontrivial finite
groups of class at most two and exponent $p$, $p$ an odd prime. Then
$G\amalg^{\mathfrak{N}_2}\!\!H$ is capable.
\label{th:anycoprod}
\end{theorem}
\begin{proof} If either $G$ or~$H$ are cyclic then the result
  follows from Lemma~\ref{lemma:coprodcyclic}. If they are both
  noncyclic, let $G$ be minimally generated by $g_1,\ldots,g_a$,
  $a>1$, and let~$H$ be minimally generated by
  $h_{a+1},\ldots,h_{a+b}$, $b>1$. Let $X_1$ be the subspace of
  $\langle v_{ji}\,|\, 1\leq i<j\leq a\rangle$ corresponding to~$G$
  under the obvious identification, and let $X_2$ be the subspace of
  $\langle v_{ji}\,|\, a+1\leq i<j\leq a+b\rangle$ corresponding
  to~$H$. Then $X_1\oplus X_2$ is the subspace of $V=\langle
  v_{ji}\,|\, 1\leq i<j\leq a+b\rangle$ that corresponds to
  $G\amalg^{{\germ N}_2} H$. The result will follow from
  Lemma~\ref{lemma:coprodgen} below.
\end{proof}

We will need some notation prior to stating
Lemma~\ref{lemma:coprodgen}. 
Let $a,b>1$. Let $V$, $W$, and
$\varphi_1,\ldots,\varphi_{n}$ be as in~(\ref{eq:defofV})--(\ref{eq:defvarphi}), corresponding to
$n=a+b$.  We decompose $V$ as $V=V_s\oplus V_m \oplus V_{\ell}$, where $V_s$
is generated by the vectors $v_{ji}$, $1\leq i<j\leq a$ (the ``small
index'' vectors); $V_{\ell}$ is generated by the vectors $v_{ji}$,
$a+1\leq i<j\leq a+b$ (the ``large index'' vectors); and $V_m$ is
generated by the vectors $v_{ji}$ with $1\leq i\leq a<j\leq b$ (the
``mixed index'' vectors). For any vector $\mathbf{v}\in V$, we write
\[ \mathbf{v} = \mathbf{v}_s + \mathbf{v}_m + \mathbf{v}_{\ell},\]
where $\mathbf{v}_s\in V_s$, $\mathbf{v}_m\in V_m$, and
$\mathbf{v}_{\ell}\in V_{\ell}$.
The idea here is that the vectors in the ``small index'' part of $V$
will correspond to the relations defining~$G$, while the vectors in
the ``large index'' part of~$V$ will correspond to relations
defining~$H$.

Similarly, we decompose $W$ as $W=W_{s} \oplus W_{m_1}\oplus
W_{m_2}\oplus W_{\ell}$; $W_s$ is generated by the basis vectors which
have all three indices smaller than or equal to $a$; $W_{m_1}$ by the
basis vectors which have exactly two indices smaller than or equal to
$a$, one larger than $a$; $W_{m_2}$ by the basis vectors which have
exactly one index smaller than or equal to $a$, and two larger
than~$a$; and $W_{\ell}$ by the basis vectors in which all three
indices are larger than $a$. As above, we can decompose any vector
$\mathbf{w}\in W$ as
\[ \mathbf{w} = \mathbf{w}_s + \mathbf{w}_{m_1} + \mathbf{w}_{m_2} +
\mathbf{w}_{\ell},\]
with $\mathbf{w}_s\in W_s$, etc.

\begin{lemma} Notation as in the previous paragraphs. Let $X_1$ be a
  subspace of $V_s$, $X_2$ a subspace of $V_{\ell}$. If $X=X_1\oplus
  X_2$, then $X=Z_X$.
\label{lemma:coprodgen}
\end{lemma}

\begin{proof} We prove this as a series of claims. 
 Let $\mathbf{w}\in W$, 
$\mathbf{w} = \mathbf{w}_s + \mathbf{w}_{m_1} + \mathbf{w}_{m_2} +
  \mathbf{w}_{\ell}$.

CLAIM 1: $\mathbf{w}\in Y_X$ if and only if each of
$\mathbf{w}_s$, $\mathbf{w}_{m_1}$, $\mathbf{w}_{m_2}$, and
$\mathbf{w}_{\ell}$ are in $Y_X$.

Indeed, note that $\mathbf{v}\in X$ if and only if $\mathbf{v}_s\in X_1$,
$\mathbf{v}_{\ell}\in X_2$, and $\mathbf{v}_m = \mathbf{0}$. Since
$\varphi_i(X_1)\subset W_s$, $\varphi_i(X_2)\subset W_{m_2}$ for
$1\leq i\leq a$; and $\varphi_j(X_1)\subset W_{m_1}$,
$\varphi_j(X_2)\subset W_{\ell}$ if $a+1\leq j\leq a+b$, the
claim follows.

CLAIM 2: $\mathbf{v}\in Z_X$ if and only if
$\mathbf{v}_s,\mathbf{v}_{\ell}\in Z_X$ and
$\mathbf{v}_m=\mathbf{0}$. 

For, let $\mathbf{v}\in
\varphi_i^{-1}(Y_X)$ for some fixed $i$,
$1\leq i\leq a+b$. Then $\mathbf{w}=\varphi_i(\mathbf{v})\in Y_X$, so
each of $\mathbf{w}_s$, $\mathbf{w}_{m_1}$, $\mathbf{w}_{m_2}$, and
$\mathbf{w}_{\ell}$ are in $Y_X$. If $i\leq a$, we have that
$\mathbf{w}_s=\varphi_i(\mathbf{v}_s)$, $\mathbf{w}_{m_1} =
\varphi_i(\mathbf{v}_m)$, $\mathbf{w}_{m_2} =
\varphi_i(\mathbf{v}_{\ell})$, and $\mathbf{w}_{\ell}=\mathbf{0}$. On
the other hand, if
$a<i\leq a+b$, then $\mathbf{w}_s = \mathbf{0}$, $\mathbf{w}_{m_1} =
\varphi_i(\mathbf{v}_s)$, $\mathbf{w}_{m_2}=\varphi_i(\mathbf{v}_m)$,
and $\mathbf{w}_{\ell} = \varphi_i(\mathbf{v}_{\ell})$. In either
case, we see that each of $\mathbf{v}_s$, $\mathbf{v}_m$, and
$\mathbf{v}_{\ell}$ lie in $\varphi_i^{-1}(Y_X)$. This proves that
$\mathbf{v}\in Z_X$ if and only if each of $\mathbf{v}_s$,
$\mathbf{v}_{\ell}$, and $\mathbf{v}_m$ are in $Z_X$. Since all
vectors in~$X$ have zero $(j,i)$ coordinate for all $1\leq i\leq
a<j\leq a+b$, so do all vectors in $Z_X$ by
Lemma~\ref{lemma:zerocoord}, which gives the condition
$\mathbf{v}_m=\mathbf{0}$.

CLAIM 3: For each $i$, $a+1\leq i\leq a+b$, 
\[ \varphi_i(V_s) \cap \Bigl\langle
\varphi_{a+1}(X_1),\ldots,\varphi_{i-1}(X_1),
\varphi_{i+1}(X_1),\ldots,\varphi_{a+b}(X_1)\Bigr\rangle 
= \bigl\{ \mathbf{0}\bigr\}.\] Indeed, for any $j$, $a+1\leq j\leq
a+b$, the nonzero coordinates in any vector of $\varphi_j(V_s)$ have
last index equal to~$j$. Thus none of the nonzero coordinates of the
vectors in the span of the $\varphi_j(X_1)$, $a+1\leq j\leq a+b$,
$j\neq i$, have a nonzero coordinate with last index equal to $i$, and
hence cannot lie in $\varphi_i(V_s)$.

CLAIM 4: $Z_X\cap V_s = X_1$. 

Let $\mathbf{v}_s\in Z_X\cap V_s$. Then
\[ \mathbf{v}_s \in \bigcap_{i=a+1}^{a+b} \varphi_i^{-1}(Y_X).\]
and thus,
\[\mathbf{v}_s \in \bigcap_{i=a+1}^{a+b}
  \varphi_i^{-1}(Y_X\cap W_{m_1}).\]
But $Y_X\cap W_{m_1} = \langle
\varphi_{a+1}(X_1),\ldots,\varphi_{a+b}(X_1)\rangle$, and from Claim~3
it follows that 
\[ \varphi_i^{-1}\left(\bigl\langle
\varphi_{a+1}(X_1),\ldots,\varphi_{b+1}(X_1)\bigr\rangle\right) =
X_1;\quad i=a+1,\ldots,a+b.\]
Therefore, $\mathbf{v}_s\in X_1$, as claimed.

CLAIM 5: For each $i$, $1\leq i\leq a$, 
\[ \varphi_i(V_{\ell}) \cap \Bigl\langle
\varphi_{1}(X_2),\ldots,\varphi_{i-1}(X_2),\varphi_{i+1}(X_2),\ldots,\varphi_a(X_2)\Bigr\rangle
= \bigl\{\mathbf{0}\bigr\}.\] 
This is analogous to Claim~3: the
nonzero coordinates in any vector in $\varphi_j(V_{\ell})$, $1\leq
j\leq a$, have middle index equal to $j$. So none of the nonzero
vectors in the span of the $\varphi_j(X_2)$, $1\leq j\leq a$, $j\neq
i$ may lie in $\varphi_i(V_{\ell})$.

CLAIM 6: $Z_X\cap V_{\ell} = X_2$. 

This follows from Claim~5 in the same way that Claim~4 follows from Claim~3.

CLAIM 7: $Z_X= X_1\oplus X_2 = X$. 

From Claims~2, 4, and~6 we have that
$Z_X = (Z_X\cap V_s)\oplus(Z_X\cap V_{\ell}) = X_1\oplus X_2 = X$,
as claimed. This proves the lemma.
\end{proof}

\subsection*{Reducing to a special case}
A more interesting result is the following:

\begin{theorem} Let $G$ be a finite noncyclic nilpotent group of class two
  and exponent~$p$, $p$ an odd prime; let $C_p$ be the cyclic group of
  order~$p$. Then
\[ G\mbox{\ is capable}\Longleftrightarrow G\oplus C_p\mbox{\ is
  capable.}\]
\label{th:cancellation}
\end{theorem}
\begin{proof} Let $G$ be minimally generated by $n$
  elements $g_1,\ldots,g_n$, $n>1$. We think of $C_p$ as generated by
  $x_{n+1}$. Let $V,W$ be the spaces corresponding to $n+1$, and let
  $X_1$ be the subspace of $V_1 = \langle v_{ji}\,|\, 1\leq i<j\leq n\rangle$
that corresponds to $G$.  The subspace that corresponds to $G\oplus
C_p$ is $X_1 \oplus \langle v_{(n+1)i}\,|\, 1\leq i\leq n\rangle$. The
statement that $G$ is capable is the statement that $X_1 = Z_1$ (where
we work only with $V_1$ and $\varphi_1,\ldots,\varphi_n$), and the
statement that $G\oplus C_p$ is capable is equivalent to saying that
$X=Z_X$ (this time working with $V$ and
$\varphi_1,\ldots,\varphi_{n+1}$). Therefore, the theorem will follow from
  Lemma~\ref{lemma:addcyclic} below.
\end{proof}

\begin{lemma} Let $n>1$, and let $V$, $W$,
  $\varphi_1,\ldots,\varphi_{n+1}$ be as in
  $(\ref{eq:defofV})$--$(\ref{eq:defvarphi})$, corresponding to
  $n+1$. Let
\begin{eqnarray*}
V_1 & = & \bigl\langle v_{ji}\,\bigm|\, 1\leq i<j\leq n\bigr\rangle,\\
W_1 & = & \bigl\langle w_{jik}\,\bigm|\, 1\leq i<j\leq n, \>\; i\leq k\leq n\bigr\rangle.
\end{eqnarray*}
Let $X_1$ be a subspace of $V_1$, and let
\begin{eqnarray*}
Y_1 &=& \Bigl\langle \varphi_1(X_1),\ldots,\varphi_n(X_1)\Bigr\rangle\\
Z_1 &=& \bigcap_{i=1}^n \varphi_i^{-1}(Y_1).
\end{eqnarray*}
Let $X$ be the subspace of~$V$ given by
\[ X = X_1 \oplus \bigl\langle v_{(n+1)i}\,|\, 1\leq i\leq n\bigr\rangle,\] 
and let $Y_X$ and $Z_X$ be given by
\begin{eqnarray*}
 Y_X & = & \Bigl\langle
 \varphi_1(X),\ldots,\varphi_n(X),\varphi_{n+1}(X)\Bigr\rangle\\
Z_X & = & \bigcap_{i=1}^{n+1} \varphi_i^{-1}(Y_X).
\end{eqnarray*}
Then $X_1=Z_1$ if and only if $X=Z$.
\label{lemma:addcyclic}
\end{lemma}

\begin{proof}
Note that all vectors $w_{(n+1)ij}$ and $w_{ji(n+1)}$, with $1\leq
i<j\leq n+1$ are in $Y_X$. For $w_{(n+1)ij}$ is the image under
  $\varphi_{j}$ of $v_{(n+1)i}$; and the image under $\varphi_i$ of
  $v_{(n+1)j}$ is $w_{(n+1)ij} - w_{ji(n+1)}$; the first summand is
already in $Y_X$, hence so is the second summand. This means that
$V_1\subset \varphi_{n+1}^{-1}(Y_X)$, since $v_{ji} =
\varphi_{n+1}^{-1}(w_{ji(n+1)})$.  Since all vectors $v_{(n+1)i}$ are
in $X$, we conclude that in fact $V = \varphi_{n+1}^{-1}(Y_X)$; (in
terms of the group-theoretic setting, all we are saying is that since
every generator commutes with $x_{n+1}$, every commutator does as well).

Let $X_2 = \langle v_{(n+1),i}\,|\, 1\leq i\leq n\rangle$.

Also note that $Y_1$ does not contain any vector which has a nonzero
coordinate with $n+1$ in any index. So we have:
\begin{eqnarray*}
Y_X & = & \bigl\langle
\varphi_1(X),\ldots,\varphi_{n+1}(X)\bigr\rangle\\
& = & \bigl\langle
\varphi_1(X_1),\ldots,\varphi_n(X_1),
\varphi_1(X_2),\ldots,\varphi_n(X_2),\varphi_{n+1}(X)\bigr\rangle\\
& = & \bigl\langle
Y_1,\varphi_1(X_2),\ldots,\varphi_n(X_2),\varphi_{n+1}(X)\big\rangle\\
& = & Y_1 \oplus \bigl\langle w_{(n+1)ij},w_{ji(n+1)}\,\bigm|\,
1\leq i<j\leq n+1\bigr\rangle.
\end{eqnarray*}

Since $Z_1$ also does not contain any vector which has a nonzero
coordinate with $n+1$ in any index, we claim that 
\begin{equation}
Z_X = Z_1 \oplus \bigl\langle v_{(n+1)i}\,\bigm|\, 1\leq i\leq
n\bigr\rangle;
\label{eq:zxaszxone}
\end{equation}
and from this, $Z_X = X \Longleftrightarrow Z_1=X_1$ will follow.

Let $\mathbf{v}\in Z_X$. We want to show that it is in
$Z_1\oplus\langle v_{(n+1)i}\rangle$. By adding the necessary
multiples of $v_{(n+1)j}$, we may assume that $\mathbf{v}\in
V_1$. Therefore, $\varphi_i(\mathbf{v})\in Y_X\cap W_1 = Y_1$ for
$i=1,2,\ldots,n$, which implies that $\mathbf{v}\in Z_1$. So the right
hand side of (\ref{eq:zxaszxone}) contains the left hand
side. Conversely, let $\mathbf{v}\in Z_1$. To show that $\mathbf{v}\in
Z_X$, we only need to show that
$\mathbf{v}\in\varphi_{n+1}^{-1}(Y_X)$.  But since $\mathbf{v}\in
Z_1\subset V_1$, this follows from our observations above.  This
proves the lemma.
\end{proof}

The interest of Theorem~\ref{th:cancellation} is that it allows us
to reduce the problem of capability to a special case that has been
considered often in the past. We do this in the following two results: 

\begin{lemma} Let $G$ be a finite nilpotent group of class two and
  exponent~$p$, with $p$ an odd prime. Then $G$ may be written as $G=
  K\oplus C_p^r$, where $K$ satisfies $Z(K)=[K,K]=[G,G]$, $r =
  \dim_{\mathbb{F}_p}(Z(G)/[G,G])$, and $C_p$ is
  the cyclic group of order~$p$.
\end{lemma}

\begin{proof}  Let $z_1,\ldots,z_r \in Z(G)$ be elements that project
onto a basis of the vector space $Z(G)/[G,G]$. Let $g_1,\ldots,g_n$ be
elements of~$G$ whose projection extends the images of
$z_1,\ldots,z_r$ to a basis of $G/[G,G]$. Let $K=\langle
g_1,\ldots,g_n\rangle$. It is well known that a set of elements of a
nilpotent group that generates the abelianization must also generate
the group, so $G$ is generated by $z_1,\ldots,z_r,g_1,\ldots,g_n$.

We have that $Z(G) = \langle z_1,\ldots,z_r\rangle \oplus [G,G]$, and
that $\langle z_1,\ldots,z_r\rangle$ is a direct summand of~$G$. By
construction, we also have $[K,K]=[G,G]$, so 
\[G= K\oplus \langle z_1,\ldots,z_r\rangle \cong K \oplus C_p^r.\]
Since $K$ is cocentral in~$G$, $Z(K)=Z(G)\cap K$, and therefore, any
element which is central in~$K$ must be a commutator in~$G$, by choice
of $z_1,\ldots,z_r$.  That is, $Z(K)=[G,G]=[K,K]$, which proves the
lemma.
\end{proof}

The following result is now immediate:

\begin{theorem} Let $G$ be a noncyclic finite nilpotent group of class two and
  exponent $p$, $p$ an odd prime. Write as $G=K\oplus C_p^r$, where
  $C_p$ is cyclic of order~$p$, $r\geq 0$, and $K$ satisfies
  $Z(K)=[K,K]$. Then $G$ is capable if and only if $K$ is
  capable.\hfill$\Box$ 
\label{th:specialcase}
\end{theorem}

One reason why Theorem~\ref{th:specialcase} is interesting is that the
condition $Z(K)=[K,K]$ is fairly strong, and has been used in the past
in the study of capability for nilpotent groups of class~$2$ and
exponent~$p$; for example, Theorem~1 in~\cite{heinnikolova}. It also
means that we may discard certain subspaces from consideration in the
linear algebra setting: any subspace that corresponds to a group in
which $Z(G)\neq [G,G]$ may be ignored. Any subspace $X$ that contains
all vectors $v_{kj}$ and $v_{ji}$, for a fixed $j$ and $1\leq
j\leq n$, $i<j<k$, for instance. 

\begin{corollary} The capability of finite groups of class two and odd
  prime exponent~$p$ is completely determined by the capability of
 finite groups~$G$ of class two and exponent~$p$ which satisfy the
 condition $Z(G)=[G,G]$, plus the observation that a nontrivial
 abelian group of exponent~$p$ is capable if and only if it is not cyclic.
\end{corollary}

\section{Some Geometry}\label{sec:geometry}

As we mentioned above, recasting our problem into a linear algebra
setting allows us to bring other tools to bear on the
problem. Specifically, in this section we give a geometric argument,
due to David McKinnon (personal communication) which implies that if a
group of exponent~$p$ and class two is ``nonabelian enough'', then it
will necessarily be capable. Although we have only been able to apply
the argument as-is to a limited number of cases, there is some hope
that similar arguments may hold more generally.  It also seems to
present a striking example of why the linear algebra setting may be
more amenable to a final solution than the group-theoretic setting by
itself.

Since the reader may not be familiar with Grassmannians, we provide a
very quick crash course on some of the main concepts. Unfortunately,
the main property we will use, stated in
Proposition~\ref{prop:finitefibers}, seems to be in that rather
awkward situation of being well-known or easy to figure out to those
who work in algebraic geometry, yet requiring some very technical
machinery to prove and not being available explicitly anywhere for
reference. I have provided only a very rough paraphrased sketch of how
one might go about proving it, and I apologize in advance for the
inconvenience.

I have taken most of the explanations that follow (sometimes verbatim)
from personal communications with Prof.~McKinnon, and I am very
grateful for his permission to include them here. I also borrow the
description of Grassmannians from~\cite{harris}.  For now, we suspend
the assumption that $V$ and~$W$ are the vector spaces defined
in~(\ref{eq:defofV}) and (\ref{eq:defofW}). We will explicitly state
when we reinstate that assumption.

\begin{definition}
Let $V$ be a vector space of dimension $n$ over a field $F$, and let
$k$ be an integer, $0\leq k\leq n$. The \textit{Grassmannian
$Gr(k,V)$} (or $Gr(k,n)$, if the field is understood from context) is
defined to be the set of all $k$-dimensional subspaces of~$V$.
\end{definition}

In particular, when $k=1$, we have that $Gr(k,V) = Gr(1,V)$ consists
of all ``lines through the origin'' in $V$. If $\dim(V)=n>0$, then this
is well known to correspond to $\mathbb{P}^{n-1}$, projective
$(n-1)$-space, by identifying $V$ with the $n$-dimensional affine
space over $F$.

Given a $k$-dimensional linear subspace $X$ of~$V$, spanned by vectors
$\mathbf{v}_1,\ldots,\mathbf{v}_k$, we can associate to $X$ the ``pure
wedge''
\[ \lambda = \mathbf{v}_1 \wedge \cdots \wedge \mathbf{v}_k \in
\bigwedge\nolimits^k(V),\] where, as usual, $\bigwedge^k(V)$ is the
$k$-th exterior power of~$V$. In this case, $X$ uniquely determines
$\lambda$ up to scalars, for if we choose a different basis then the
corresponding vector would simply be $\mathbf{v}_1\wedge\cdots\wedge
\mathbf{v}_k$ multiplied by the determinant of the change of basis
matrix. This gives a well-defined map of sets from $Gr(k,V)$ to the
projective space $\mathbb{P}(\bigwedge^k(V))$. It can be shown that
this map is an inclusion, called the Pl\"ucker embedding, which
identifies $Gr(k,V)$ with a subvariety of projective space; this
endows $Gr(k,V)$ with the structure of an algebraic variety.

The property of Grassmannians that we will use is the following:

\begin{prop} Let $\germ{V}$ be any variety (over the field~$F$), and
let $f\colon Gr(k,V)\to\germ{V}$ be a regular map. If ${\rm
char}(F)\neq 2$, then $f$ is either constant, or has finite fibers.
\label{prop:finitefibers}
\end{prop}

\begin{remark} A regular map is a map of algebraic varieties that is
  defined everywhere, and which is locally (in the Zariski topology)
  represented by quotients of polynomials.
\end{remark}

\noindent\textit{Sketch of proof.} The proof of
  Proposition~\ref{prop:finitefibers} proceeds in two steps. First,
  one proves that the Picard group of a Grassmannian is infinite
  cyclic, and then one uses a standard argument to show that an
  infinite cyclic Picard group implies finite fibers for regular maps
  in any variety. For the first part, one may invoke the Theorem in
  pp.~32 of~\cite{mumford} to obtain a description of the homogeneous
  coordinate ring of $Gr(k,V)$ and deduce from it that the Picard
  group is infinite cyclic, or else argue by covering $Gr(k,V)$ with
  large open affine sets (whose Picard group are shown to be trivial by
  examining the corresponding Chow groups) and then invoking Proposition~6.5
  in~\cite{hartshorne} to show that the Picard group must be a
  quotient of the infinite cyclic group. Since $Gr(k,V)$ is a
  projective variety, it is known that the Picard group cannot be
  torsion, and therefore it must be infinite cyclic. For the second
  part, a nonconstant regular map $f\colon Gr(k,V)\to {\germ V}$ from
  the Grassmannian to an arbitrary variety ${\germ V}$ induces a map
  $f^*\colon{\rm Pic}({\germ V})\to{\rm Pic}(Gr(k,V))$ given by
  $f^*(D) = f^{-1}(D)$. Since ${\rm Pic}(Gr(k,V))$ is infinite cyclic,
  then a generator must be ample by Theorem~7.6 in~\cite{hartshorne},
  and so by Prop.~7.5, also from~\cite{hartshorne}, it follows that
  for every very ample $D\in {\rm Pic}({\germ V})$, either $f^*(D)$ or
  $-f^*(D)$ is ample. As a special case of Kleiman's Criterion
  (Prop.~1.27(a) in~\cite{debarre}) the support of an ample divisor
  intersects every curve in the variety, and the support of $f^*(D)$ is
  the same as the support of $-f^*(D)$. Thus, for every very ample
  $D\in{\rm Pic}({\germ V})$, the support of $f^*(D)$ intersects any
  curve on $Gr(k,V)$. Now assume that $f\colon Gr(k,V)\to{\germ V}$ is
  a regular map which is not constant, and let $y$ be any point in the
  image of $f$; one picks a very ample divisor $D$ which intersects
  the image but does not contain~$y$ (which is always possible), and
  then notes that the support of $f^*(D)=f^{-1}(D)$ is disjoint from
  $f^{-1}(y)$; so $f^{-1}(y)$ cannot contain any curves, and so must
  be $0$-dimensional (i.e. a finite union of points). Therefore $f$
  has finite fibers.\hfill$\Box$

\smallskip

Assume that we have $n$ linear transformations
$\varphi_1,\ldots,\varphi_n\colon V\to W$. For every element $X\in
Gr(k,V)$, there is a subspace $Y_X$ of $W$ given by
\[ Y_X = \Bigl\langle \varphi_1(X),\ldots,\varphi_n(X)\Bigr\rangle.\]
Fix $X$, let $\{\mathbf{v}_1,\ldots,\mathbf{v}_k\}$ be a basis for $X$, and let $m={\rm
dim}(Y_X)$. Then we can identify $Y_X$ as a point in $Gr(m,W)$ by
setting it equal to
\begin{equation}
\sum_I \bigwedge_{(i,j)\in I} \varphi_i(\mathbf{v}_j)
\label{eq:bigsum}
\end{equation}
where $I$ ranges over all subsets of
$\{1,\ldots,n\}\times\{1,\ldots,k\}$ of cardinality $m$. By
construction, each pure wedge in the sum is a scalar multiple of a
unique pure wedge associated to $Y_X$, and at least one of the
summands is nonzero; by avoiding certain degenerate choices of bases
for $X$, we can ensure that the sum itself is nonzero, and thus yields
a nonzero pure wedge which is an element of $Gr(m,W)$.

Suppose further that there is a neighbourhood $U$ of $X$ in $Gr(k,V)$
(in the Zariski topology) such that for all $X'\in U$, the
corresponding subspace $Y_{X'}$ has dimension $m=\dim(Y_X)$.  Then by
choosing bases for elements of $U$ in a suitably well-behaved manner
(explicitly, by choosing, near $X$, a Zariski-local basis of sections of
the tautological bundle on $Gr(k,V)$), we see that the correspondence
$X\mapsto Y_X$ is in fact a rational map from $Gr(k,V)$ to $Gr(m,W)$.
That is, $X\mapsto Y_X$ is defined on a dense open subset of $X$ (in
the Zariski topology), and is locally defined by quotients of
polynomials, obtained by expanding the wedge product in the definition
above.  If $\dim(Y_X)=m$ for all $X$ of dimension $k$, then we can
take $U=Gr(k,V)$, and hence our rational map is defined everywhere,
and is therefore a regular map from $Gr(k,V)$ to $Gr(m,W)$.

\begin{remark} The discussion above is a bit more general than we will actually need
  for the limited cases in which we have been able to apply
  Theorem~\ref{th:mckinnon} below. In those cases, we have that
  $\dim(Y_X) = nk$ for all $X$ of a given dimension~$k$, and
  therefore, one need only choose a specific basis
  $\mathbf{v}_1,\ldots,\mathbf{v}_k$ of~$X$, and define the morphism
  $Gr(k,V)\to Gr(nk,W)$ by mapping $\mathbf{v}_1\wedge\cdots\wedge
  \mathbf{v}_k$ to
\[ \varphi_1(\mathbf{v}_1)\wedge\cdots\wedge
  \varphi_1(\mathbf{v}_k)\wedge \varphi_2(\mathbf{v}_1) \wedge
  \varphi_2(\mathbf{v}_2)\wedge\cdots\wedge \varphi_n(\mathbf{v}_k)
  \in Gr(nk,W).\]
Since each component of this wedge is determined by a linear
  transformation, this will be a regular map.
\end{remark}

With these notions in hand, we now return to the situation we have
associated to the group theoretic setting. Once again, we assume that
$V$ and~$W$ are defined by (\ref{eq:defofV}) and~(\ref{eq:defofW}),
respectively.

\begin{theorem}[David McKinnon~\cite{persdave}]
Let $n\geq 2$, and let $V$, $W$, and $\varphi_1,\ldots,\varphi_n$ be 
as in~$(\ref{eq:defofV})$--$(\ref{eq:defvarphi})$. Let $\overline{V} =
V\otimes_{\mathbb{F}_p}\overline{\mathbb{F}_p}$ and $\overline{W}=
W\otimes_{\mathbb{F}_p}\overline{\mathbb{F}_p}$, where
$\overline{\mathbb{F}_p}$ is the algebraic closure of
$\mathbb{F}_p$. Extend $\varphi_i$ to maps from $\overline{V}$ to
$\overline{W}$ in the obvious way.  Let $k,m$ be integers,
$0<k<\dim(\overline{V})$, $0<m<\dim(\overline{W})$. If $Y_X\in
Gr(m,\overline{W})$ whenever $X\in Gr(k,\overline{V})$, then $X=Z_X$
for all $X\in Gr(k,V)$.
\label{th:mckinnon}
\end{theorem}

\begin{proof}
By using the correspondence sketched above, since every subspace $X$
of $\overline{V}$ of dimension~$k$ corresponds to a subspace of
$\overline{W}$ of dimension~$m$, we can define a regular map
from $Gr(k,\overline{V})$ to $Gr(m,\overline{W})$.

Since our map is clearly not constant, it must have finite fibers by
Proposition~\ref{prop:finitefibers}.  We claim that in fact it is
one-to-one. For, assume to the contrary, that there are two distinct
subspaces $X_1,X_2$ of~$\overline{V}$, both of dimension~$k$, such
that $Y_1=Y_2$ in~$\overline{W}$. Then for every subspace $X'$
contained in~$\langle X_1,X_2\rangle$, we will necessarily have that
$\varphi_r(X')\subset Y_1$. In particular, for any subspace $X'$
contained in~$\langle X_1,X_2\rangle$, if ${\rm dim}(X') = k$ then
$Y_{X'}=Y_1$. However, there are infinitely many such $X'$, which would
mean that $Y_1$ has an infinite fiber. As this is impossible, we
conclude that the map $Gr(k,\overline{V}) \to Gr(m,\overline{W})$ is
one-to-one. Therefore, the corresponding map $Gr(k,V)\to Gr(m,W)$ is
also one-to-one.

Let $X$ be a subspace of~$V$ of dimension~$k$, and consider $Z_X$,
defined as in~(\ref{eq:ZsubX}).
If $X'$ is any subspace of
dimension~$k$ contained in~$Z_X$, then $Y_{X'} = Y_X$. Since the
correspondence from $Gr(k,V)$ to $Gr(m,W)$ is
one-to-one, there can be only one subspace of dimension~$k$ contained
in $Z_X$. Therefore, $X=Z_X$, as claimed.
\end{proof}

When do we have the condition given in Theorem~\ref{th:mckinnon}? The
following result gives a partial answer to that question (and,
unfortunately, it is ``seldom''):

\begin{prop} Let $n>2$, and let $\overline{V}$ and $\overline{W}$ be as
 in Theorem~\ref{th:mckinnon}. Let $X$ be a subspace of $\overline{V}$.
\begin{itemize}
\item[(i)] If ${\rm dim}(X)=1$, then ${\rm dim}(Y_X) = n$.
\item[(ii)] If ${\rm dim}(X) = 2$, then ${\rm dim}(Y_X) = 2n$.
\item[(iii)] If $2<k<n$, then there exist subspaces $X_1$ and~$X_2$ of
  $\overline{V}$ (in fact, subspaces that correspond to subspaces of
  $V$), with ${\rm dim}(X_1)={\rm dim}(X_2)=k$, but ${\rm
  dim}(Y_1)\neq {\rm dim}(Y_2)$.
\end{itemize}
\label{prop:limits}
\end{prop}

\begin{proof} For (iii), let
\begin{eqnarray*}
X_1&=&\langle  v_{21},\ldots,v_{(k+1)1}\rangle\\
X_2&=&\langle v_{21}, v_{31},v_{32}, v_{51},\ldots,v_{(k+1)1}\rangle.
\end{eqnarray*}
Then it is easy to verify that:
\begin{eqnarray*}
Y_1 &=& \langle w_{r1s}\,|\, 1\leq s\leq n, 2\leq r\leq (k+1)\rangle\\
Y_2 &=& \langle w_{r1s}, w_{32u}\,|\, 1\leq s\leq n,
r=2,3,5,\ldots,(k+1), 2\leq u\leq n \rangle.
\end{eqnarray*}
so ${\rm dim}(Y_1) = kn$ and ${\rm dim}(Y_2) = kn-1$.

For (i), assume that $X=\langle\mathbf{v}\rangle$,
$\mathbf{v}\neq\mathbf{0}$. If the dimension is not equal to $n$, then
the set $\{\varphi_1(\mathbf{v}),\ldots,\varphi_n(\mathbf{v})\}$ is
linearly dependent. Let $i_0$ be the first index such that
$\varphi_{i_0}(\mathbf{v})$ is a linear combination of the previous
vectors. Note that $i_0>1$.

Neither $r$ nor $s$ can equal $i_0$, because only
$\varphi_{i_0}(\mathbf{v})$ would have a nonzero coordinate 
indexed by a triple that includes two $i_0$'s. 

If $s<i_0<r$, then $\varphi_{i_0}(\mathbf{v})$ has nonzero
$(r,s,i_0)$ coordinate. Since the only other way to get a vector with
nonzero $(r,s,i_0)$ coordinate is by applying $\varphi_s$ to a
vector with nonzero $(r,i_0)$ coordinate, then $\mathbf{v}$ must have
nonzero $(r,i_0)$ coordinate; but we have already noted this is impossible.

If $i_0<s<r$, then $\varphi_{i_0}(\mathbf{v})$ has nonzero
$(r,i_0,s)$ and $(s,i_0,r)$ coordinates; so to express it as a linear
combination of other images of $\mathbf{v}$, we must be
using the image under $\varphi_r$. Since $r>i_0$, this is impossible
by choice of $i_0$.

Thus, we conclude that $r<i_0$, in which case
$\varphi_{i_0}(\mathbf{v})$ involves a nonzero $(r,s,i_0)$
coordinate. As above, that means that $\varphi_s(\mathbf{v})$ must
have nonzero $(r,s,i_0)$ coordinate, which means that $\mathbf{v}$
has nonzero $(i_0,r)$ coordinate. But this is impossible as well, as
we already noted.
Therefore, $Y_X$ must be of dimension $n$ when ${\rm dim}(X)=1$, as claimed.

For (ii), we proceed as above. Let $X=\langle
\mathbf{v}_1,\mathbf{v}_2\rangle$, of dimension~$2$. Assume that the
set 
\[\{
\varphi_1(\mathbf{v}_1),\varphi_1(\mathbf{v}_2),\varphi_2(\mathbf{v}_1),
\varphi_2(\mathbf{v}_2),\ldots, \varphi_n(\mathbf{v}_1),\varphi_n(\mathbf{v}_2)\}\]
is linearly dependent. By exchanging $\mathbf{v}_1$ and
$\mathbf{v}_2$, and replacing $\mathbf{v}_1$ with a linear combination
of $\mathbf{v}_1$ and $\mathbf{v}_2$ if necessary, we may assume that the first
vector in the set which is a linear combination of the previous ones
is $\varphi_{i_0}(\mathbf{v}_1)$ for some $i_0>1$. 

Let $(r,s)$ be a nonzero coordinate of $\mathbf{v}_1$.  As above, we
cannot have $i_0\in\{r,s\}$, and we must have $s<r<i_0$.
Then $\varphi_{i_0}(\mathbf{v}_1)$ has nonzero $(r,s,i_0)$
coordinate. To obtain that, we must be using the image under
$\varphi_s$ of vector with nonzero $(i_0,r)$ coordinate; that vector
cannot be $\mathbf{v}_1$, so it is $\mathbf{v}_2$ that has a nonzero
$(i_0,r)$ coordinate. Since $\varphi_s(\mathbf{v}_2)$ has nonzero
$(i_0,s,r)$ coordinate as well, and $\varphi_{i_0}(\mathbf{v}_1)$ does
not, we must also be using the image under $\varphi_r$ of a vector
with nonzero $(i_0,s)$ coordinate; once again, this cannot be
$\mathbf{v}_1$, so it is $\mathbf{v}_2$ which has a nonzero $(i_0,s)$
coordinate as well. But then $\varphi_s(\mathbf{v}_2)$ has a nonzero
$(i_0,s,s)$ coordinate, which is not the case for
$\varphi_{i_0}(\mathbf{v}_1)$ and which cannot be cancelled with any
other available vector. This is a contradiction. Thus we conclude that
${\rm dim}(Y_X)=2n$ when ${\rm dim}(X)=2$.
\end{proof}

Note that if $G$ is a finite nilpotent group of class~$2$ and
exponent~$p$, then $G$ is minimally $n$-generated where
$n=\dim_{\mathbb{F}_p}(G^{\rm ab})$. For certainly we cannot generate
$G$ with fewer elements, and any generating set for $G^{\rm ab}$ lifts
to a generating set for~$G$.

\begin{corollary} Let $G$ be a finite nonabelian nilpotent group of
class 2 and exponent~$p$ an odd prime, minimally generated by
$x_1,\ldots,x_n$, $n>2$. If
\[ \dim_{\mathbb{F}_p}\bigl([G,G]\bigr) \geq \binom{n}{2} - 2\]
then $G$ is capable.
\label{cor:nonabenough}
\end{corollary}

\begin{proof} Since $\#[G,G] = p^{\binom{n}{2}}/\#N$, the condition
  guarantees that $\#N\leq p^2$, so if $X$ corresponds to $G$, then
  ${\rm dim}(X)\leq 2$. The result now follows from
  Theorem~\ref{th:mckinnon}, Proposition~\ref{prop:limits}, and the
  trivial case of $X=\{\mathbf{0}\}$.
\end{proof}

\begin{remark} Although Proposition~\ref{prop:limits} shows that the
  applicability of Theorem~\ref{th:mckinnon}, as is, is limited, it may be
  possible to use similar ideas to obtain other results. It seems
  likely that ``most'' subspaces $X$ will have ${\rm dim}(Y_X)$ equal
  to a fixed number, except for some degenerate exceptions where the
  dimension is smaller. If a map could be defined from a sufficiently
  nice subvariety of ${\rm Gr}(k,\overline{V})$, the conclusion may
  still follow.
\end{remark}

Combining Corollary~\ref{cor:nonabenough} with
Theorem~\ref{th:specialcase} yields:

\begin{theorem} Let $G$ be a finite noncyclic group of class two and
  exponent~$p$, and let $k=\dim_{\mathbb{F}_p}(G/Z(G))$. If 
\[ \dim_{\mathbb{F}_p}\bigl([G,G]\bigr) \geq \binom{k}{2} - 2\]
then $G$ is capable.
\label{th:generalsufficient}
\end{theorem}

\begin{proof} We may write $G = K \oplus C_p^r$, where
  $r=\dim_{\mathbb{F}_p}(Z(G)/[G,G])$, and $K$ satisfies
  $Z(K)=[K,K]=[G,G]$. Then $G$ is capable if and only if $K$ is
  capable, and we apply Corollary~\ref{cor:nonabenough} to $K$, noting
  that $\dim_{\mathbb{F}_p}(K^{\rm ab}) =
  \dim_{\mathbb{F}_p}(G/Z(G))=k$. 
\end{proof}

The result is strong enough to settle the 3-generated case:

\begin{theorem} Let $G$ be a $3$-generated group of class 2 and
  exponent~$p$. Then $G$ is either cyclic or capable.
\end{theorem}

\begin{proof} Assume that $G$ is noncyclic and $3$-generated. If $G$
  is abelian, then it is capable by Theorem~\ref{th:trivialcases}. So
  we may assume that $G$ is not abelian. If $G$ is minimally
  $2$-generated, then it is the nonabelian group of order $p^3$, which
  is isomorphic to the $2$-nilpotent product of two cyclic groups of
  order~$p$, so we again conclude it is capable by
  Theorem~\ref{th:trivialcases}. The only remaining case has $G$
  minimally $3$-generated, and $\#[G,G]\geq p$, which is capable by
  Corollary~\ref{cor:nonabenough}. 
\end{proof}

Note that we already have an example of a minimally 4-generated nilpotent
group of class~$2$ and exponent~$p$ which is not capable.

\section{Final comments}\label{sec:fincomments}

As another application of Theorem~\ref{th:specialcase}, we use a
theorem of Heineken and Nikolova to give a bound on the dimension of
the subpace $X$ that we need to consider.

\begin{theorem}[Heineken and Nikolova, Theorem~1
in~\cite{heinnikolova}] Assume that $G$ is capable, of
exponent~$p$, and $Z(G)=[G,G]$. If $Z(G)$ is of rank~$k$, then the
rank of $G^{\rm ab}$ is at most $2k+\binom{k}{2}$.
\label{th:heinnikolova}
\end{theorem}

First, a very easy application of Theorem~\ref{th:specialcase} allows
us to drop the hypothesis $Z(G)=[G,G]$, provided we replace $Z(G)$
with $[G,G]$ and $G^{\rm ab}$ with $G/Z(G)$:  

\begin{corollary}
Assume that $G$ is capable, of exponent~$p$ and class at most~$2$. If
$[G,G]$ is of rank~$k$, then the rank of $G/Z(G)$ is at most
$2k+\binom{k}{2}$. 
\label{cor:extensionhn}
\end{corollary}

\begin{proof} Write $G=K\oplus C_p^r$, with $C_p$ cyclic of order~$p$
and $K$ is a subgroup satisfying $Z(K)=[K,K]$. Since $G$ is capable, $K$ is
capable. Thus, Theorem~\ref{th:heinnikolova} applies to~$K$. Now
simply note that $K^{\rm ab}\cong G/Z(G)$ and that $Z(K)=[K,K]=[G,G]$.
\end{proof}

By turning the theorem ``upside down'' as it were, we can give a lower
bound on the size of $[G,G]$ in terms of the size of a minimal
generating set for $G$.

\begin{theorem} Let $G$ be a finite nilpotent group of class
at most~$2$ and exponent an odd prime~$p$. If $G$ is capable and
$\dim_{\mathbb{F}_p}(G/Z(G)) = n$, then 
\[ \dim_{\mathbb{F}_p}([G,G]) \geq \left\lceil
\frac{-3+\sqrt{9+8n}}{2}\,\right\rceil,\] 
where $\lceil x\rceil$ denotes the smallest integer greater than or
equal to~$x$.
\label{th:generalnecessary}
\end{theorem}

\begin{proof} From Corollary~\ref{cor:extensionhn}, we know that if
$\dim_{\mathbb{F}_p}([G,G])=k$, then we will have that
$\dim_{\mathbb{F}_p}(G/Z(G))=n\leq 2k+\binom{k}{2}$. We 
transform the inequality $n\leq 2k+\binom{k}{2}$ into an inequality
for~$k$ in terms of~$n$; the inequality is equivalent to
$k^2 + 3k - 2n \geq 0$; and since both
$k$ and $n$ must be nonnegative integers, this gives that
\[ k \geq \left\lceil\frac{-3+\sqrt{9+8n}}{2}\,\right\rceil,\]
as claimed.
\end{proof}

If we restrict to groups~$G$ satisfying $Z(G)=[G,G]$, noting that
$\dim_{\mathbb{F}_p}\bigl([G,G]\bigr) = \binom{n}{2} -
\dim_{\mathbb{F}_p}(X)$ (where~$X$ is determined by~$G$), we obtain:

\begin{corollary} The capability of finite groups of class two and odd
prime exponent~$p$ is completely determined by considering $V$, $W$,
$\varphi_1,\ldots,\varphi_n$ as in
$(\ref{eq:defofV})$--$(\ref{eq:defvarphi})$, and subspaces $X$ of~$V$
with
\[ \dim_{\mathbb{F}_p}(X) \leq \binom{n}{2} -
\left\lceil\frac{-3+\sqrt{9+8n}}{2}\,\right\rceil \]
and which do not correspond to groups~$G$ for which $Z(G)\neq [G,G]$.
\end{corollary}

In view of these results, and particularly of
Theorem~\ref{th:anycoprod}, it seems that capability for groups of
class two depends on there not being too many relations among the
commutators (that is, the subspace $X$ not being ``too big'').
We have not quite succeeded in closing the gap between the necessary
condition in Theorem~\ref{th:generalnecessary} and the sufficient
condition in Theorem~\ref{th:generalsufficient}. But, echoing the
comments in~\cite{baconkappe}, it seems reasonable to hope that this
gap may be closed soon.

\section*{Aknowledgments}
It is my very great pleasure to thank David McKinnon for providing
most of the geometry in Section~\ref{sec:geometry}, and specifically
Theorem~\ref{th:mckinnon} and the details behind my sketchy paraphrase
of the second half of the proof of
Proposition~\ref{prop:finitefibers}. The first part of that sketch comes from details provided by 
N.I.~Shepherd-Barron and by
Mike Roth, and their help is very much
appreciated.

\end{document}